# Evaluation of Ramanujan Continued Fractions

**Nikos Bagis**
Stenimahou 5 Edessa Pellas
58200 Greece
bagkis@hotmail.com

**Abstract**
In this paper we present experimental ways of evaluating Ramanujan`s quantities which as someone can see are related with algebraic numbers. The good thing with algebraic numbers is that can be found in a closed form, from there approximations, using Mathematica. In this way we produce new formulas and give new ideas for to prove new theorems.

We use the notation

$$b_0 + \cfrac{a_1}{b_1 + \cfrac{a_2}{b_2 + \cfrac{a_3}{b_3 + \ldots}}} := b_0 + \frac{a_1}{b_1 +}\frac{a_2}{b_2 +}\frac{a_3}{b_3 +}\ldots$$

thus for example

$$1 + \frac{1}{1+}\frac{1}{1+}\frac{1}{1+}\frac{1}{1+}\ldots = \frac{\sqrt{5}+1}{2}$$

For $|q|<1$, we define the "Rogers-Ramanujan continued fraction" by

$$R(q) := \frac{q^{1/5}}{1}+\frac{q}{1+}\frac{q^2}{1+}\frac{q^3}{1+}\ldots$$

and

$$R^*(q) := \frac{1}{1+}\frac{q}{1+}\frac{q^2}{1+}\frac{q^3}{1+}\ldots$$

For $|q|<1$, we denote the "Rogers-Ramanujan functions,"

$$(a;q)_n := \prod_{k=0}^{n-1}(1-aq^k) \tag{1}$$

Define the Ramanujan`s function

$$f(-q) := (q;q)_\infty$$

For $|q|<1$, Rogers and Ramanujan show that

$$q^{-1/5}R(q) = R^*(q) := \frac{(q;q^5)_\infty (q^4;q^5)_\infty}{(q^2;q^5)_\infty (q^3;q^5)_\infty} = \prod_{n=1}^\infty (1-q^n)^{X_2(n)}$$

where

$$X_2(n) = \begin{cases} 1, & \text{if } n = 1 \bmod 5 \\ -1, & \text{if } n = 2 \bmod 5 \\ -1, & \text{if } n = 3 \bmod 5 \\ 1, & \text{if } n = 4 \bmod 5 \\ 0, & \text{if } n = 0 \bmod 5 \end{cases} = \left(\frac{n}{5}\right)$$

Also there exists following relations which we shall use to recover some values of the Rogers Ramanujan continued fraction.

With $f$ defined by (2),

$$\frac{1}{R(q)} - 1 - R(q) = \frac{f(-q^{1/5})}{q^{1/5} f(-q^5)} \tag{2}$$

$$\frac{1}{R^5(q)} - 11 - R^5(q) = \frac{f^6(-q)}{q \cdot f^6(-q^5)} \tag{3}$$

The following identity is valid

$$\log(R^*(x)) = -\sum_{n=1}^\infty \frac{x^n}{n} \sum_{d|n} X_2(d) d, \quad |x| < 1$$

The above relation shows clearly that the derivatives in zero of the function $\log(R^*(q))$ are

$$-\Gamma(n) \sum_{d|n} X_2(d) d$$

It's also holds

$$R(e^{-x}) = e^{-x/5} \frac{\vartheta_4(3ix/4, e^{-5x/2})}{\vartheta_4(ix/4, e^{-5x/2})}, \quad x > 0$$

$$R'(q) = \frac{2 \cdot 2^{1/3} k^{1/3} (k')^{4/3} K^2}{5\pi^2 q} R(q) \cdot \sqrt[6]{\frac{1}{R(q)^5} - 11 - R(q)^5}$$

For the $K, k, k', \vartheta_4$, see [W,W]:

$$R(e^{-2\pi}) = -\frac{1}{2} - \frac{\sqrt{5}}{2} + \sqrt{\frac{5+\sqrt{5}}{2}}$$

and

$$R'(e^{-2\pi}) = 8\sqrt{\frac{2}{5}\left(9 + 5\sqrt{5} - 2\sqrt{50 + 22\sqrt{5}}\right)} \frac{e^{2\pi}}{\pi^3} \Gamma\left(\frac{5}{4}\right)^4$$

$$R(e^{-x}) = \exp\left(-\frac{x}{5} - \sum_{n=1}^{\infty} \frac{1}{n} \frac{e^{4nx} - e^{3nx} - e^{2nx} + e^{nx}}{e^{5nx} - 1}\right), \quad x > 0$$

Also the following formula is valid:

$$\frac{1}{1+} \frac{-e^{-x}}{1+e^{-x}} + \frac{-e^{-3x}}{1+e^{-2x}} + \frac{-e^{-5x}}{1+e^{-3x}} + \frac{-e^{-7x}}{1+e^{-4x}} + \ldots = \exp\left(-\sum_{n=1}^{\infty} \frac{e^{-nx}}{n} \sum_{d|n} Y_2(d) d\right),$$

where $Y_2(n) = \begin{cases} 1, & \text{if } n = 1 \bmod 3 \\ -1, & \text{if } n = 2 \bmod 3 \\ 0, & \text{if } n = 0 \bmod 3 \end{cases}$

Thus as above:

$$\frac{1}{1+} \frac{-e^{-x}}{1+e^{-x}} + \frac{-e^{-3x}}{1+e^{-2x}} + \frac{-e^{-5x}}{1+e^{-3x}} + \frac{-e^{-7x}}{1+e^{-4x}} + \ldots = \exp\left(-\sum_{n=1}^{\infty} \frac{1}{n} \frac{e^{nx} - e^{2nx}}{e^{3nx} - 1}\right)$$

For the continued fraction

$$H(x) = \frac{e^{-x/2}}{(1+e^{-x})+} \frac{e^{-2x}}{(1+e^{-3x})+} \frac{e^{-4x}}{(1+e^{-5x})+} \frac{e^{-6x}}{(1+e^{-7x})+} \frac{e^{-8x}}{(1+e^{-9x})+} \ldots =$$

$$= \exp\left(-x/2 - \sum_{n=1}^{\infty} \frac{1}{n} \frac{e^{7nx} - e^{5nx} - e^{3nx} + e^{nx}}{e^{8nx} - 1}\right)$$

exists the following relation

If $ab = \pi^2$, then

$$\left(1 + \sqrt{2} + H(a)\right)\left(1 + \sqrt{2} + H(b)\right) = 2(2 + \sqrt{2}) \tag{4}$$

We also have $H(\pi\sqrt{r}) = $ Algebraic, where $r$ positive rational, for example we have

$$H\left(\frac{\pi}{2}\right) = \sqrt{1 + 2\sqrt{2} - 2\sqrt{2 + \sqrt{2}}}$$

$$H\left(\frac{\pi\sqrt{2}}{2}\right) = \sqrt{3 + 2\sqrt{2} - 2\sqrt{4 + 3\sqrt{2}}}$$

-A good thing with relative series and products is that they give us a nice precision to calculate continued fractions with Mathematica. They also give us the related polynomials and closed forms of the continued fractions and prove new theorems, set new problems as in relation (4).

Another interesting formula is:

When $q = e^{-\pi\sqrt{r}}$ $a, p, r \in \mathbb{Q}_+^*$

$$q^{\frac{2a}{3} - 4p + \frac{4p^2}{a}} (q^{a-p}; q^a)_\infty^8 (q^p; q^a)_\infty^8 = \text{Algebraic} \qquad (5)$$

This relation explains most of the algebraic properties of Ramanujan`s fractions.

Also another interesting result is

$$-\int_a^b f(-q)^4 q^{-5/6} dq =$$

$$= \left[ 6R(x)^{5/6} \sum_{n=0}^{\infty} \frac{\left(\frac{1}{6}\right)_n^2}{\left(\frac{7}{6}\right)_n n!} \left(\frac{11 + 5\sqrt{5}}{2}\right)^n {}_2F_1\left[\begin{array}{c}\frac{1}{6}, -n \\ \frac{5}{6} - n\end{array}; \frac{11 - 5\sqrt{5}}{11 + 5\sqrt{5}}\right] R(x)^{5n} \right]_a^b \qquad (6)$$

and

$$2\pi \int_a^b \eta(i\tau)^4 d\tau =$$

$$-6R(e^{2\pi y})^{5/6} F_1\left[\tfrac{1}{6}, \tfrac{1}{6}, \tfrac{1}{6}, \tfrac{7}{6}, \tfrac{1}{2}(11 - 5\sqrt{5})R(e^{-2\pi y})^5, \tfrac{1}{2}(11 + 5\sqrt{5})R(e^{-2\pi y})^5\right]_a^b \qquad (7)$$

Where $F_1[a, b_1, b_2, c, x, y] = \sum_{m=0}^{\infty} \sum_{n=0}^{\infty} \frac{(a)_{m+n}(b_1)_m(b_2)_n}{m! n! (c)_{m+n}} x^m y^n$ is Appell function

M.L.Glasser also prove that

$$\int_0^1 f(-q)^4 q^{-5/6} dq = \pi 2^{1/6} (\sqrt{5} - 1)^{5/6} {}_2F_1\left[\begin{array}{c}1/6, 1/6 \\ 1\end{array}; \frac{1}{2}(-123 + 55\sqrt{5})\right] \qquad (8)$$

$$\int_0^1 f(-q^5)^4 q^{-1/6} dq = \frac{1}{8 \cdot 2^{1/6}} \pi (\sqrt{5} - 1)^{25/6} {}_2F_1\left[\begin{array}{c}5/6, 5/6 \\ 1\end{array}; \frac{1}{2}(-123 + 55\sqrt{5})\right] \qquad (9)$$

From (13) and (14) we get

$$\int_0^1 \eta(ix)^4 dx = 1/2 \left(\frac{\sqrt{5} - 1}{2}\right)^{5/6} {}_2F_1\left[\begin{array}{c}1/6, 1/6 \\ 1\end{array}; \frac{1}{2}(-123 + 55\sqrt{5})\right] \qquad (10)$$

where $\eta$ denotes the Dedekind Eta function see [W,W].

$${}_2F_1\left[\begin{array}{c}5/6, 5/6 \\ 1\end{array}; \frac{1}{2}(-123 + 55\sqrt{5})\right] = 1/5 \left(\frac{\sqrt{5} + 1}{2}\right)^{10/3} {}_2F_1\left[\begin{array}{c}1/6, 1/6 \\ 1\end{array}; \frac{-123 + 55\sqrt{5}}{2}\right]$$

-Another result is

$$R(q) = \frac{1}{2}\left(\sqrt{x^2+2x+5}-x-1\right), \quad x = \frac{\eta\left(-\frac{\log(q)}{10\pi}\right)}{\eta\left(-\frac{i5\log(q)}{2\pi}\right)},$$

From which we get

$$\frac{1}{\pi}\frac{d}{d\tau}\log(x_1(\tau)) = \frac{4\eta(i\tau)^4\sqrt{x_1(\tau)^2+2x_1(\tau)+5}}{x_1(\tau)F\left(\sqrt{x_1(\tau)^2+2x_1(\tau)+5}-x_1(\tau)-1\right)} \quad (11)$$

where $F(w) = \dfrac{10}{(-11+32/w^5-w^5/32)^{1/6}}$, $x_1(\tau) = \dfrac{\eta(\tau i/5)}{\eta(5\tau i)}$

**Theta functions evaluation**

If $n,a,b,c,r \in \mathbb{Q}_+^*$, $q_n = e^{-\pi n\sqrt{r}}$ then

$$\left(q_n^{\frac{b^2-4ac}{4a}}\sqrt{\frac{\pi}{K(k_r)}}\sum_{v=-\infty}^{\infty}q_n^{av^2+bv+c}\right)^8 = \text{Algebraic} \quad (12)$$

For example let $q = e^{-\pi}$, $K(k_1) = \dfrac{\Gamma(1/4)^2}{4\sqrt{\pi}}$ then

$$\left(q^{1/24}\sqrt{\frac{\pi}{K(k_1)}}\sum_{k=-\infty}^{\infty}q^{(3k-1)k/2}\right)^8 = \frac{14}{27}+\frac{8}{3\sqrt{3}}+\frac{1}{2}\sqrt{\frac{2048}{243}+\frac{3584}{243\sqrt{3}}}$$

**Other Continued Fractions**

When $q = e^{\pi\sqrt{r}}$, $r \in \mathbb{Q}_+^*$ we have (see [L,W]):

**i)**

$$\frac{q^{1/3}}{1+}\frac{q+q^2}{1+}\frac{q^2+q^4}{1+}\frac{q^3+q^6}{1+}\ldots = q^{1/3}\frac{(q;q^6)_\infty(q^5,q^6)_\infty}{(q^3;q^6)_\infty^2} = \text{Algebraic}$$

**ii)** The same thing holds also for the Rogers Ramanujan Continued Fraction.

$$\frac{q^{1/5}}{1+}\frac{q}{1+}\frac{q^2}{1+}\frac{q^3}{1+}\ldots = q^{1/5}\frac{(q;q^5)_\infty(q^4,q^5)_\infty}{(q^2;q^5)_\infty(q^3;q^5)} = \text{Algebraic}$$

and

**iii)**

$$\frac{q^{1/2}}{(1+q)+}\frac{q^2}{(1+q^3)+}\frac{q^4}{(1+q^5)+}\frac{q^6}{(1+q^7)+}\ldots = q^{1/2}\frac{(q;q^8)_\infty(q^7,q^8)_\infty}{(q^3;q^8)_\infty(q^5;q^8)_\infty} = \text{Algebraic}$$

$$(-q^2;q^2)_\infty^8 = q^{-2/3}\text{Algebraic}$$

and

$$(-q;q^2)_\infty^8 = q^{1/3}\text{Algebraic}$$

Hence

$$q\left(\frac{(-q^2;q^2)_\infty}{(-q;q^2)_\infty}\right)^8 = q\left(\frac{1}{1+}\frac{q}{1+}\frac{q^2+q}{1+}\frac{q^3}{1+}\frac{q^4+q^2}{1+}\frac{q^5}{1+}\ldots\right)^8 = \text{Algebraic}$$

**iv)**

$$q^{1/3}\frac{(q;q^2)_\infty}{\{(q^3;q^6)_\infty\}^3} = \frac{q^{1/3}}{1+}\frac{q+q^2}{1+}\frac{q^2+q^4}{1+}\frac{q^3+q^6}{1+}\ldots = \text{Algebraic}$$

**v)** We have

$$\sum_{k=0}^\infty (-c)^k q^{k(k+1)/2} = \frac{1}{1+}\frac{cq}{1+}\frac{c(q^2-q)}{1+}\frac{cq^3}{1+}\frac{c(q^4-q^2)}{1+}\ldots$$

we set

$$M(c,q) := \frac{1}{1-}\frac{cq}{1+}\frac{c(q-q^2)}{1-}\frac{cq^3}{1+}\frac{c(q^2-q^4)}{1-}\ldots = \sum_{k=0}^\infty c^k q^{k(k+1)/2}$$

then

$$M(c,q) + \frac{1}{c}\cdot M(1/c,q) = \sum_{k=-\infty}^\infty c^k q^{k(k+1)/2}$$

If $q = e^{-\pi n\sqrt{r}}$, $c = -q^a$, $a \in \mathbb{Q}_+^*$ then

$$M(q^a,q) + q^{-a}\cdot M(q^{-a},q) = \sum_{k=-\infty}^\infty q^{k^2/2+(a+1/2)k}$$

Hence

$$M(q^a,q) + q^{-a}\cdot M(q^{-a},q) = \sqrt{\frac{K}{\pi}}q^{-\left(\frac{a^2}{2}+\frac{a}{2}+\frac{1}{8}\right)}\text{Algebraic}$$

Or better

$$\frac{1}{1-}\frac{q^{a+1}}{1+}\frac{q^{a+1}(1-q)}{1-}\frac{q^{a+3}}{1+}\frac{q^{a+2}(1-q^2)}{1-}\ldots$$

$$+\frac{q^{-a}}{1-}\frac{q^{-a+1}}{1+}\frac{q^{-a+1}(1-q)}{1-}\frac{q^{-a+3}}{1+}\frac{q^{-a+2}(1-q^2)}{1-}+\ldots = \sqrt{\frac{K}{\pi}}q^{-\left(\frac{a^2}{2}+\frac{a}{2}+\frac{1}{8}\right)}\text{Algebraic}$$

For $a = 0$ we get

$$\frac{1}{1-}\frac{q}{1+}\frac{(q-q^2)}{1-}\frac{q^3}{1+}\frac{(q^2-q^4)}{1-}\ldots = \sum_{k=-\infty}^{\infty} q^{k(k+1)/2} = q^{-1/8}\sqrt{\frac{K(k_r)}{\pi}} \cdot \text{Algebraic}$$

An "exact" evaluation of $M(c,q)$ is when $c=-q^a$, $a$ is odd integer (the proof is easy):

$$\frac{q^{\frac{(a+1)^2}{4}}}{1-}\frac{q^{a+2}}{1-}\frac{q^a(q^4-q^2)}{1-}\frac{q^{a+6}}{1-}\frac{q^a(q^8-q^4)}{1-}\frac{q^{a+10}}{1-}\frac{q^a(q^{12}-q^6)}{1-}\ldots = \frac{1}{2} - \sum_{k=0}^{(a-1)/2} q^{k^2} + \frac{\vartheta_3(q)}{2}$$

where $\vartheta_3(q) = \sum_{k=-\infty}^{\infty} q^{k^2}$ and

**vi)**

$$\frac{q^{1/2}}{1-q} + \frac{q(1-q)^2}{(1-q)(q^2+1)} + \frac{q(1-q^3)^2}{(1-q)(q^4+1)} + \frac{q(1-q^5)^2}{(1-q)(q^6+1)} + \ldots =$$

$$= q^{1/2}\frac{(q^4,q^4)_\infty^2}{(q^2;q^4)_\infty^2} = \frac{K}{\pi} \cdot \text{Algebraic}$$

**The Rogers Ramanujan type continued fractions and their first order derivatives**

$$R_1(q) = \frac{q^{1/5}}{1+}\frac{q}{1+}\frac{q^2}{1+}\frac{q^3}{1+}\ldots = q^{1/5}\frac{(q;q^5)_\infty(q^4,q^5)_\infty}{(q^2;q^5)_\infty(q^3;q^5)} = \text{Alg}$$

$$R_2(q) = \frac{q^{1/3}}{1+}\frac{q+q^2}{1+}\frac{q^2+q^4}{1+}\frac{q^3+q^6}{1+}\ldots = q^{1/3}\frac{(q;q^6)_\infty(q^5,q^6)_\infty}{(q^3;q^6)_\infty^2} = \text{Alg}$$

$$R_3(q) = \frac{q^{1/2}}{(1+q)+}\frac{q^2}{(1+q^3)+}\frac{q^4}{(1+q^5)+}\frac{q^6}{(1+q^7)+}\ldots = q^{1/2}\frac{(q;q^8)_\infty(q^7,q^8)_\infty}{(q^3;q^8)_\infty(q^5;q^8)_\infty} = \text{Alg}$$

have all derivative such as

$$R_{1,2,3}'(q)\frac{q\cdot\pi^2}{K^2} = \text{Algebraic}$$

Thus we observe that

$$\frac{d}{dQ}\left(Q^{-\frac{(p_1-p_2)}{2}+\frac{(p_1^2-p_2^2)}{2a}}\frac{(Q^{a-p_1};Q^a)_\infty(Q^{p_1};Q^a)_\infty}{(Q^{a-p_2};Q^a)_\infty(Q^{p_2};Q^a)_\infty}\right)_{Q=q=e^{-\pi\sqrt{r}}} = \frac{K^2}{q\pi^2}\text{Algebraic}$$

also

$$\frac{d}{dQ}\left(Q^{\frac{a}{12}-\frac{p}{2}+\frac{p^2}{2a}}(Q^{a-p};Q^a)_\infty(Q^p;Q^a)_\infty\right)_{Q=q=e^{-\pi\sqrt{r}}} = \frac{K^2}{q\pi^2}\text{Algebraic}$$

For example

$$\frac{d}{dQ}\left(Q^{1/8}\frac{(Q^3;Q^4)_\infty(Q;Q^4)_\infty}{(Q^2;Q^4)_\infty^2}\right)_{Q=q=e^{-\pi}} = \frac{e^\pi \Gamma(1/4)^4}{64\cdot 2^{5/8}\pi^3}$$

$$\frac{d}{dQ}\left(Q^{1/5}\frac{(Q^4;Q^5)_\infty(Q;Q^5)_\infty}{(Q^3;Q^5)_\infty(Q^2;Q^5)_\infty}\right)_{Q=q=e^{-\pi}} = \frac{e^\pi \Gamma(1/4)^4}{16\cdot \pi^3}\rho$$

Where $\rho$ is root of

$$16 - 240t^2 + 800t^3 - 2900t^4 - 6000t^5 - 6500t^6 + 17500t^7 + 625t^8 = 0$$

and

$$\frac{d}{dQ}\left(Q^{-1/24}(Q^3;Q^4)_\infty(Q;Q^4)_\infty\right)_{Q=q=e^{-\pi}} = \frac{-e^\pi \Gamma(1/4)^4}{32\cdot 2^{7/8}\pi^3}$$

### References.